\documentclass[12pt,english,reqno]{amsart}
\usepackage{a4wide}
\usepackage{amsmath,amssymb,amsfonts,amsthm}

\usepackage[
  pdfauthor={Bernd C. Kellner},
  pdfkeywords={Asymptotic constants, Glaisher-Kinkelin constant,
  Bernoulli number, factorials, Gamma function, Riemann zeta function},
  pdftitle={On Asymptotic Constants Related to Products of Bernoulli
  Numbers and Factorials},
  colorlinks,bookmarksnumbered]{hyperref}

\newcommand{\pdiv}{\mid}
\newcommand{\notdiv}{\nmid}
\newcommand{\mods}[1]{({\rm mod\,} #1)}
\newcommand{\divsum}{{}\!'\,}

\DeclareMathOperator{\real}{Re}
\DeclareMathOperator{\imag}{Im}
\DeclareMathOperator*{\spn}{span}

\newcommand{\BN}{\widehat{B}}
\newcommand{\RR}{\mathbb{R}}
\newcommand{\CC}{\mathbb{C}}
\newcommand{\NN}{\mathbb{N}}
\newcommand{\MA}{\mathcal{A}}
\newcommand{\MB}{\mathcal{B}}
\newcommand{\MC}{\mathcal{C}}
\newcommand{\MF}{\mathcal{F}}

\numberwithin{equation}{section}

\newtheorem{prop}{Proposition}[section]
\newtheorem{theorem}[prop]{Theorem}
\newtheorem{corl}[prop]{Corollary}
\newtheorem{lemma}[prop]{Lemma}

\theoremstyle{definition}
\newtheorem{defin}[prop]{Definition}
\newtheorem{result}[prop]{Result}

\theoremstyle{remark}
\newtheorem{remark}[prop]{Remark}

\title[On Asymptotic Constants ...]
{On Asymptotic Constants Related to Products of Bernoulli Numbers and Factorials}
\author{Bernd C. Kellner}
\date{}
\address{Mathematisches Institut, Universit\"at G\"ottingen, Bunsenstr.~3–-5,
37073 G\"ottingen, Germany}
\email{bk@bernoulli.org}

\begin{document}

\begin{abstract}
We discuss the asymptotic expansions of certain products of Bernoulli numbers
and factorials, e.g.,
\[
   \prod_{\nu=1}^n |B_{2\nu}| \quad \text{and} \quad \prod_{\nu=1}^n (k \nu)!^{\,\nu^r}
     \quad \text{as} \quad n \to \infty
\]
for integers $k \geq 1$ and $r \geq 0$. Our main interest is to determine exact
expressions, in terms of known constants, for the asymptotic constants of these
expansions and to show some relations among them.
\end{abstract}

\maketitle

\section{Introduction}

Let $B_n$ be the $n$th Bernoulli number. These numbers are defined by
\[
   \frac{z}{e^z-1} = \sum_{n=0}^\infty B_n \frac{z^n}{n!},
     \qquad |z| < 2 \pi,
\]
where $B_n = 0$ for odd $n > 1$. The Riemann zeta function $\zeta(s)$
is defined by
\begin{equation} \label{eqn-def-zeta}
  \zeta(s) = \sum_{\nu=1}^\infty \nu^{-s}
    = \prod_{p} (1-p^{-s})^{-1}, \quad s \in \CC, \,\, \real s > 1.
\end{equation}
By Euler's formula we have for even positive integers $n$ that
\begin{equation} \label{eqn-euler-zeta}
   \zeta(n) = - \frac12 \frac{(2\pi i)^n}{n!} \, B_n.
\end{equation}
\smallskip

Products of Bernoulli numbers occur in certain contexts in number theory.
For example, the Minkowski--Siegel mass formula states for positive integers
$n$ with $8 \pdiv n$ that
\[
   M(n) = \frac{|B_k|}{2k} \,\, \prod_{\nu=1}^{k-1} \frac{|B_{2\nu}|}{4\nu},
     \quad n=2k,
\]
which describes the mass of the genus of even unimodular positive definite
$n \times n$ matrices, for details see \cite[p.~252]{Koecher98}.
We introduce the following constants which we shall need further on.

\begin{lemma} \label{lem-const-c1}
There exist the constants
\begin{align*}
  \MC_1 &= \prod_{\nu = 2}^\infty \zeta(\nu)    = 2.2948565916... \,, \\
  \MC_2 &= \prod_{\nu = 1}^\infty \zeta(2\nu)   = 1.8210174514... \,, \\
  \MC_3 &= \prod_{\nu = 1}^\infty \zeta(2\nu+1) = 1.2602057107... \,.
\end{align*}
\end{lemma}

\begin{proof}
We have $\log (1+x) < x$ for real $x>0$. Then
\begin{equation} \label{eqn-estimate-zeta}
   \log \prod_{\nu = 1}^\infty \zeta(2\nu)
     = \sum _{\nu = 1}^\infty \log \zeta(2\nu)
     < \sum _{\nu = 1}^\infty ( \zeta(2\nu)-1) = \frac34 .
\end{equation}
The last sum of \eqref{eqn-estimate-zeta} is well known and follows by
rearranging in geometric series, since we have absolute convergence.
We then obtain that $\pi^2/6 < \MC_2 < e^{3/4}$, $\zeta(3) < \MC_3 < \MC_2$,
and $\MC_1 = \MC_2 \MC_3$.
\end{proof}

To compute the infinite products above within a given precision, one can use
the following arguments. A standard estimate for the partial sum of $\zeta(s)$
is given by
\[
   \zeta(s) - \sum_{\nu=1}^N \nu^{-s} < \frac{N^{1-s}}{s-1},
     \quad s \in \RR, \,\, s > 1.
\]
This follows by comparing the sum of $\nu^{-s}$ and the integral of $x^{-s}$ in
the interval $(N,\infty)$. Now, one can estimate the number $N$ depending on
$s$ and the needed precision. However, we use a computer algebra system, that
computes $\zeta(s)$ to a given precision with already accelerated built-in
algorithms. Since $\zeta(s) \to 1$ monotonically as $s \to \infty$, we next
have to determine a finite product that suffices the precision. From above, we
obtain
\begin{equation} \label{eqn-estimate-zeta-2}
   \zeta(s) - 1 < 2^{-s} \left(1 + \frac{2}{s-1}\right),
     \quad s \in \RR, \,\, s > 1.
\end{equation}
According to \eqref{eqn-estimate-zeta} and \eqref{eqn-estimate-zeta-2},
we then get an estimate for the remainder of the infinite product by
\[
   \log \prod_{\nu > N'} \zeta(\nu) < 2^{-N'+\varepsilon}
\]
where we can take $\varepsilon = 3/N'$; the choice of $\varepsilon$ follows by
$2^x \geq 1 + x \log 2$ and \eqref{eqn-estimate-zeta-2}.
\medskip

We give the following example where the constant $\MC_1$ plays an important
role; see Finch \cite{Finch03}. Let $a(n)$ be the number of non-isomorphic
abelian groups of order $n$. The constant $\MC_1$ equals the average of the
numbers $a(n)$ by taking the limit. Thus, we have
\[
   \MC_1 = \lim_{N \to \infty} \frac{1}{N} \sum_{n=1}^N a(n).
\]

By definition the constant $\MC_2$ is connected with values of the Riemann zeta
function on the positive real axis. Moreover, this constant is also connected
with values of the Dedekind eta function
\[
   \eta( \tau ) = e^{\pi i \tau/12} \prod_{\nu = 1}^\infty
     ( 1 - e^{2\pi i \nu \tau} ), \quad \tau \in \CC, \,\, \imag \tau > 0,
\]
on the upper imaginary axis.

\begin{lemma}
The constant $\MC_2$ is given by
\[
   1 / \MC_2 = \prod_p p^{\frac{1}{12}} \,\, \eta \!
     \left( i \frac{\log p}{\pi} \right)
\]
where the product runs over all primes.
\end{lemma}

\begin{proof}
By Lemma \ref{lem-const-c1} and the Euler product \eqref{eqn-def-zeta} of
$\zeta(s)$, we obtain
\[
   \MC_2 = \prod_{\nu = 1}^\infty \prod_{p} (1-p^{-2\nu})^{-1}
     = \prod_{p} \prod_{\nu = 1}^\infty (1-p^{-2\nu})^{-1}
\]
where we can change the order of the products because of absolute convergence.
Rewriting $p^{-2\nu} = e^{2\pi i\nu\tau}$ with $\tau = i \log p \, / \pi$
yields the result.
\end{proof}

We used \textsc{Mathematica} \cite{Wolfram} to compute all numerical values in
this paper. The values were checked again by increasing the needed precision
to 10 more digits.

\section{Preliminaries}

We use the notation $f \sim g$ for real-valued functions when
$\lim_{x \to \infty} f(x)/g(x) = 1$. As usual, $O(\cdot)$ denotes Landau's
symbol. We write $\log f$ for $\log(f(x))$.

\begin{defin}
Define the linear function spaces
\[
   \Omega_n = \spn\limits_{0 \leq \nu \leq n}
     \{ x^\nu, \, x^\nu \log x \}, \quad n \geq 0,
\]
over $\RR$ where $f \in \Omega_n$ is a function $f: \RR^+ \to \RR$.
Let
\[
   \Omega_\infty = \bigcup_{n \geq 0} \Omega_n.
\]
Define the linear map $\psi : \Omega_\infty \to \RR$ which gives the constant
term of any $f \in \Omega_\infty$. For the class of functions
\begin{equation} \label{eqn-omega}
   F(x) = f(x) + O(x^{-\delta}), \quad f \in \Omega_n,
     \quad n \geq 0, \quad \delta > 0,
\end{equation}
define the linear operator $[\cdot]: C(\RR^+;\RR) \to \Omega_\infty$ such that
$[F]=f$ and $[F] \in \Omega_n$. Then $\psi([F])$ is defined to be the
asymptotic constant of $F$.
\end{defin}

We shall examine functions $h: \, \NN \to \RR$ which grow exponentially; in
particular these functions are represented by certain products. Our problem is
to find an asymptotic function $\tilde{h}: \RR^+ \to \RR$ where
$h \sim \tilde{h}$. If $F = \log \tilde{h}$ satisfies \eqref{eqn-omega}, then
we have $[\log \tilde{h}] \in \Omega_n$ for a suitable $n$ and we identify
$[\log \tilde{h}] = [\log h] \in \Omega_n$ in that case.

\begin{lemma} \label{lem-fg-omega}
Let $f \in \Omega_n$ where
\[
   f(x) = \sum_{\nu=0}^n ( \alpha_\nu \, x^\nu + \beta_\nu \, x^\nu \log x )
\]
with coefficients $\alpha_\nu, \beta_\nu \in \RR$.
Let $g(x)=f(\lambda x)$ with a fixed $\lambda \in \RR^+$.
Then $g \in \Omega_n$ and $\psi(g) = \psi(f) + \beta_0 \log \lambda$.
\end{lemma}

\begin{proof}
Since $g(x)=f(\lambda x)$ we obtain
\[
   g(x) = \sum_{\nu=0}^n ( \alpha_\nu \, (\lambda x)^\nu
     + \beta_\nu \, (\lambda x)^\nu ( \log \lambda + \log x ) ).
\]
This shows that $g \in \Omega_n$. The constant terms are $\alpha_0$ and
$\beta_0 \log \lambda$, thus $\psi(g) = \psi(f) + \beta_0 \log \lambda$.
\end{proof}

\begin{defin}
For a function $f: \, \RR^+ \to \RR$ we introduce the notation
\[
   f(x) = \sum_{\nu \geq 1} \divsum f_\nu(x)
\]
with functions $f_\nu: \, \RR^+ \to \RR$ in case $f$ has a divergent series
expansion such that
\[
   f(x) = \sum_{\nu=1}^{m-1} f_\nu(x) + \theta_m(x) f_m(x),
     \quad \theta_m(x) \in (0,1), \quad m \geq N_f,
\]
where $N_f$ is a suitable constant depending on $f$.
\end{defin}

Next we need some well known facts which we state without proof,
cf.\ \cite{Graham94}.

\begin{prop}
Let
\[
   H_0 = 0, \qquad H_n = \sum_{\nu=1}^n \frac{1}{\nu}, \quad n \geq 1,
\]
be the $n$th harmonic number. These numbers satisfy
$H_n = \gamma + \log n + O (n^{-1})$ for $n \geq 1$,
where $\gamma = 0.5772156649...$ is Euler's constant.
\end{prop}

\begin{prop}[Stirling's series] \label{prop-gamma-stirling}
The Gamma function $\Gamma(x)$ has the divergent series expansion
\[
   \log \Gamma(x+1) = \frac12 \log(2\pi) + \left( x+\frac12 \right) \log x
     - x + \sum_{\nu \geq 1} \divsum \frac{B_{2\nu}}{2\nu (2\nu-1)} \,
     x^{-(2\nu-1)}, \ x > 0.
\]
\end{prop}

\begin{remark}
When evaluating the divergent series given above, we have to choose a suitable
index $m$ such that
\[
   \sum_{\nu \geq 1} \divsum \frac{B_{2\nu}}{2\nu (2\nu-1)} \, x^{-(2\nu-1)}
     = \sum_{\nu = 1}^{m-1} \frac{B_{2\nu}}{2\nu (2\nu-1)} \, x^{-(2\nu-1)}
     + \theta_m(x) R_m(x)
\]
and the remainder $|\theta_m(x) R_m(x)|$ is as small as possible. Since
$\theta_m(x) \in (0,1)$ is not effectively computable in general, we have to
use $|R_m(x)|$ instead as an error bound. Sch\"afke and Finsterer
\cite{Schaefke90}, among others, showed that the so-called Lindel\"of error
bound $L = 1$ for the estimate $L \geq \theta_m(x)$ is best possible for
positive real $x$.
\end{remark}

\begin{prop} \label{prop-fac-gamma}
If $\alpha \in \RR$ with $0 \leq \alpha < 1$, then
\[
   \prod_{\nu=1}^n (\nu - \alpha) = \frac{\Gamma(n+1-\alpha)}{\Gamma(1-\alpha)}
     \,\,\sim\,\, \frac{\sqrt{2\pi}}{\Gamma(1-\alpha)} \, \left( \frac{n}{e}
     \right)^n \, n^{\frac12-\alpha} \quad \text{as} \quad n \to \infty.
\]
\end{prop}

\begin{prop}[Euler] \label{prop-euler-gamma}
Let $\Gamma(x)$ be the Gamma function. Then
\[
   \prod_{\nu=1}^{n-1} \Gamma \! \left( \frac{\nu}{n} \right) =
     \frac{(2\pi)^{\frac{n-1}{2}}}{\sqrt{n}}.
\]
\end{prop}

\begin{prop}[Glaisher \cite{Glaisher87}, Kinkelin \cite{Kinkelin60}]
\label{prop-glaisher}
Asymptotically, we have
\[
   \prod_{\nu=1}^n \nu^\nu \,\,\sim\,\, \MA \,\,
     n^{\frac12 n(n+1)+\frac{1}{12}} \, e^{-\frac{n^2}{4}}
     \quad \text{as} \quad n \to \infty
\]
where $\MA = 1.2824271291...$ is the Glaisher--Kinkelin constant,
which is given by
\[
   \log \MA = \frac{1}{12} - \zeta'(-1)
     = \frac{\gamma}{12} + \frac{1}{12} \log (2\pi) - \frac{\zeta'(2)}{2\pi^2}.
\]
\end{prop}

Numerous digits of the decimal expansion of the Glaisher--Kinkelin constant
$\MA$ are recorded as sequence A074962 in OEIS \cite{Sloane}.

\section{Products of factorials}

In this section we consider products of factorials and determine their
asymptotic expansions and constants. For these asymptotic constants we
derive a divergent series representation as well as a closed formula.

\begin{theorem} \label{thm-asympt-fk}
Let $k$ be a positive integer. Asymptotically, we have
\[
   \prod_{\nu=1}^n (k \nu)!  \,\sim\,
     \MF_k \, \MA^k \, (2\pi)^{\frac14} \!
     \left( \frac{k\,n}{e^{3/2}} \right)^{\frac{k}{2} n(n+1)} \!\!
     \left( 2\pi k e^{k/2-1} \, n \right)^{\frac{n}{2}}
     n^{\frac14+\frac{k}{12}+\frac{1}{12k}}
     \quad \! \text{as} \quad \! n \to \infty
\]
with certain constants $\MF_k$ which satisfy
\[
   \log \MF_k = \frac{\gamma}{12k} +
     \sum_{j \geq 2} \divsum \frac{B_{2j} \,
     \zeta(2j-1)}{2j (2j-1) \, k^{2j-1}}.
\]
Moreover, the constants have the asymptotic behavior that
\[
   \lim_{k \to \infty} \MF_k = 1, \quad
     \lim_{k \to \infty} \MF_k^k = e^{\gamma/12},
     \quad \text{and} \quad
     \prod_{k=1}^n  \MF_k \,\,\sim\,\,
     \MF_\infty \,\, n^{\gamma/12}
     \quad \text{as} \quad n \to \infty
\]
with
\[
   \log \MF_\infty = \frac{\gamma^2}{12} +
     \sum_{j\geq2} \divsum \frac{B_{2j} \, \zeta(2j-1)^2}{2j (2j-1)}.
\]
\end{theorem}

\begin{theorem} \label{thm-formula-fk}
If $k$ is a positive integer, then
\[
   \log \MF_k = -\left( k+\frac{1}{k} \right) \log \MA + \frac{1}{12k}
     - \frac{1}{12k} \log k + \frac{k}{4} \log (2\pi)
     - \sum_{\nu=1}^{k-1} \frac{\nu}{k} \log \Gamma \!
       \left( \frac{\nu}{k} \right).
\]
\end{theorem}

We will prove Theorem \ref{thm-formula-fk} later, since we shall need
several preparations.

\begin{proof}[Proof of Theorem~\ref{thm-asympt-fk}]
Let $k \geq 1$ be fixed. By Stirling's approximation,
see Proposition \ref{prop-gamma-stirling}, we have
\begin{equation} \label{eqn-log-kn-fac}
   \log (k\nu)! = \frac12 \log(2\pi) + \left( k \nu + \frac12 \right)
     \log(k\nu) - k\nu + f(k\nu)
\end{equation}
where we can write the remaining divergent sum as
\[
   f(k\nu) = \frac{1}{12k\nu} + \sum_{j \geq 2} \divsum
     \frac{B_{2j}}{2j (2j-1) \, (k\nu)^{2j-1}}.
\]
Define $S(n)=1+\cdots+n=n(n+1)/2$. By summation we obtain
\begin{align*}
   \sum_{\nu=1}^n \log (k\nu)! = \frac{n}{2} \log(2\pi k)
     &+ \frac12 \log n! - k S(n) + k S(n) \log k \\
     &+ \, k \sum_{\nu=1}^n \nu \log \nu + \sum_{\nu=1}^n f(k\nu).
\end{align*}
The term $\frac12 \log n!$ is evaluated again by \eqref{eqn-log-kn-fac}.
Proposition \ref{prop-glaisher} provides that
\[
   k \sum_{\nu=1}^n \nu \log \nu = k \log \MA + k S(n) \log n
     + \frac{k}{12} \log n - \frac{k}{2} \left( S(n)-\frac{n}{2} \right)
     + O(n^{-\delta})
\]
with some $\delta > 0$. Since $\lim_{n \to \infty} H_n - \log n = \gamma$,
we asymptotically obtain for the remaining sum that
\begin{equation} \label{eqn-log-fk-1}
   \lim_{n \to \infty} \left( \sum_{\nu=1}^n f(k\nu)
       - \frac{1}{12k} \log n \right)
     = \frac{\gamma}{12k} + \sum_{j \geq 2} \divsum
       \frac{B_{2j} \, \zeta(2j-1)}{2j (2j-1) \, k^{2j-1}}
     =: \log \MF_k.
\end{equation}
Here we have used the following arguments. We choose a fixed index $m > 2$
for the remainder of the divergent sum. Then
\begin{equation} \label{eqn-log-fk-2}
   \lim_{n \to \infty} \sum_{\nu=1}^n \theta_m(k\nu)
     \frac{B_{2m}}{2m (2m-1) \, (k\nu)^{2m-1}}
   = \eta_m \frac{B_{2m} \, \zeta(2m-1)}{2m (2m-1) \, k^{2m-1}}
\end{equation}
with some $\eta_m \in (0,1)$, since $\theta_m(k\nu) \in (0,1)$ for all
$\nu \geq 1$. Thus, we can write \eqref{eqn-log-fk-1} as an asymptotic series
again. Collecting all terms, we finally get the asymptotic formula
\begin{align*}
   \sum_{\nu=1}^n \log (k\nu)! = \log \MF_k &+ k \log \MA
        + \frac14 \log(2\pi) + \, k S(n) \left( -\frac32 + \log (k n) \right) \\
     &+ \frac{n}{2} \left( \log(2\pi k) + \frac{k}{2} - 1 + \log n \right) \\
     &+ \left( \frac14 + \frac{k}{12} + \frac{1}{12k} \right) \log n
        + O(n^{-\delta'})
\end{align*}
with some $\delta' > 0$. Note that the exact value of $\delta'$ does not
play a role here. Now, let $k$ be an arbitrary positive integer.
From \eqref{eqn-log-fk-1} we deduce that
\begin{equation} \label{eqn-log-k-fk}
   \log \MF_k = \frac{\gamma}{12k} + O(k^{-3})
     \qquad \text{and} \qquad
     k \log \MF_k = \frac{\gamma}{12} + O(k^{-2}).
\end{equation}
The summation of \eqref{eqn-log-fk-1} yields
\begin{equation} \label{eqn-sum-log-fk-1}
   \sum_{k=1}^n \log \MF_k = \frac{\gamma}{12} H_n
     + \sum_{k=1}^n \sum_{j \geq 2} \divsum
     \frac{B_{2j} \, \zeta(2j-1)}{2j (2j-1) \, k^{2j-1}}.
\end{equation}
Similar to \eqref{eqn-log-fk-1} and \eqref{eqn-log-fk-2}, we can write again:
\begin{equation} \label{eqn-sum-log-fk-2}
   \lim_{n \to \infty} \left( \sum_{k=1}^n \log \MF_k - \frac{\gamma}{12}
       \log n \right)
     = \frac{\gamma^2}{12} + \sum_{j\geq2} \divsum \frac{B_{2j} \,
       \zeta(2j-1)^2}{2j (2j-1)}
     =: \log \MF_\infty. \qedhere
\end{equation}
\end{proof}

The case $k=1$ of Theorem \ref{thm-asympt-fk} is related to the so-called
Barnes $G$-function, cf.~\cite{Adamchik01}. Now we shall determine exact
expressions for the constants $\MF_k$. For $k \geq 2$ this is more complicated.

\begin{lemma} \label{lem-const-f1}
We have $\MF_1 = (2\pi)^{\frac14} \, e^{\frac{1}{12}} / \MA^2$.
\end{lemma}

\begin{proof}
Writing down the product of $n!$ repeatedly in $n+1$ rows, one observes by
counting in rows and columns that
\begin{equation} \label{eqn-prod-fac-f1}
   n!^{n+1} = \prod_{\nu=1}^n \nu! \, \prod_{\nu=1}^n \nu^\nu.
\end{equation}
From Proposition \ref{prop-gamma-stirling} we have
\[
   (n+1) \log n! = \frac{n+1}{2} \log (2\pi) - n(n+1) +
      (n+1)\left( n + \frac12 \right) \log n + \frac{1}{12} + O(n^{-1}).
\]
Comparing the asymptotic constants of both sides of \eqref{eqn-prod-fac-f1}
when $n \to \infty$, we obtain
\[
   (2\pi)^{\frac12} \, e^{\frac{1}{12}}
     = \MF_1 \, \MA \, (2\pi)^{\frac14} \cdot \MA
\]
where the right side follows by Theorem \ref{thm-asympt-fk} and
Proposition \ref{prop-glaisher}.
\end{proof}

\begin{prop} \label{prop-fprod-kl}
Let $k, l$ be integers with $k \geq 1$. Define
\[
   F_{k,l}(n) := \prod_{\nu=1}^n (k\nu-l)!
     \quad \text{for} \quad 0 \leq l < k.
\]
Then $[ \log F_{k,l} ] \in \Omega_2$ and
$F_{k,0}(n) \cdots F_{k,k-1}(n) = F_{1,0}(kn)$. Moreover
\[
   F_{k,l}(n) / F_{k,l+1}(n) = k^n \prod_{\nu=1}^n
     \left( \nu-\frac{l}{k} \right)
     \quad \text{for} \quad 0 \leq l < k-1
\]
and $[ \log(F_{k,l}/ F_{k,l+1}) ] = [ \log F_{k,l} ] - [ \log F_{k,l+1} ]
\in \Omega_1$.
\end{prop}

\begin{proof}
We deduce the proposed products from $(k\nu-l)!/(k\nu-(l+1))! = k\nu-l$ and
\begin{equation} \label{eqn-prod-Fkl}
   \prod_{\nu=1}^n (k\nu)!(k\nu-1)!\cdots(k\nu-(k-1))!
     = \prod_{\nu=1}^{kn} \nu!.
\end{equation}
Proposition \ref{prop-fac-gamma} shows that $[ \log(F_{k,l}/ F_{k,l+1}) ] \in
\Omega_1$. Since the operator $[\cdot]$ is linear, it follows that
\begin{equation} \label{eqn-log-Fkl}
   [ \log(F_{k,l}/ F_{k,l+1}) ] = [ \log F_{k,l} - \log F_{k,l+1} ]
     = [ \log F_{k,l} ] - [ \log F_{k,l+1} ] \in \Omega_1.
\end{equation}
From Theorem \ref{thm-asympt-fk} we have $[ \log F_{k,0} ] \in \Omega_2$.
By induction on $l$ and using \eqref{eqn-log-Fkl}
we derive that $[ \log F_{k,l} ] \in \Omega_2$ for $0 < l < k$.
\end{proof}

\begin{lemma} \label{lem-matrix-mk}
Let $k$ be an integer with $k \geq 2$. Define the $k \times k$ matrix
\[
   M_k := \left( \begin{array}{cccccc}
     1 & -1     \\
       & 1 & -1 \\
       &   & \ddots & \ddots \\
       &   &        & 1 & -1 \\
     1 & 1 & \cdots & 1 & 1  \\
   \end{array} \right)
\]
where all other entries are zero. Then $\det M_k = k$ and the matrix inverse
is given by $M_k^{-1} = \frac{1}{k} \widetilde{M_k}$ with
\[
   \widetilde{M_k} = \left( \begin{array}{ccccccc}
     k-1 & k-2 & k-3 & \cdots & 2 & 1 & 1 \\
      -1 & k-2 & k-3 & \cdots & 2 & 1 & 1 \\
      -1 &  -2 & k-3 & \cdots & 2 & 1 & 1 \\
      \vdots & \vdots & \vdots & & \vdots & \vdots & \vdots \\
      -1 &  -2 &  -3 & \cdots & 2 & 1 & 1 \\
      -1 &  -2 &  -3 & \cdots & -(k-2) & 1 & 1 \\
      -1 &  -2 &  -3 & \cdots & -(k-2) & -(k-1) & 1 \\
   \end{array} \right).
\]
\end{lemma}

\begin{proof}
We have $\det M_2 = 2$. Let $k \geq 3$. We recursively deduce
by the Laplacian determinant expansion by minors on the first column that
\[
   \det M_k = (-1)^{1+1} \det M_{k-1} + (-1)^{1+k} \det T_{k-1}
\]
where the latter matrix $T_{k-1}$ is a lower triangular matrix having $-1$
in its diagonal. Therefore
\[
   \det M_k = \det M_{k-1} + (-1)^{1+k} \cdot (-1)^{k-1} = k-1 + 1 = k
\]
by induction on $k$. Let $I_k$ be the $k \times k$ identity matrix.
The equation $M_k \cdot \widetilde{M_k} = k \, I_k$ is easily verified
by direct calculation, since $M_k$ has a simple form.
\end{proof}

\begin{proof}[Proof of Theorem~\ref{thm-formula-fk}]
The case $k=1$ agrees with Lemma \ref{lem-const-f1}. For now, let $k \geq 2$.
We use the relations between the functions $F_{k,l}$, resp.\ $\log F_{k,l}$,
given in Proposition \ref{prop-fprod-kl}. Since $[\log F_{k,l}] \in \Omega_2$,
we can work in $\Omega_2$. The matrix $M_k$ defined in Lemma \ref{lem-matrix-mk}
mainly describes the relations given in \eqref{eqn-prod-Fkl} and
\eqref{eqn-log-Fkl}. Furthermore we can reduce our equations to $\RR$ by
applying the linear map $\psi$, since we are only interested in the asymptotic
constants. We obtain the linear system of equations
\[
   M_k \cdot x = b \,, \qquad x, b \in \RR^k,
\]
where
\[
   x = \left( \psi([ \log F_{k,0} ]), \ldots,
       \psi([ \log F_{k,k-1} ]) \right)^T
\]
and $b=(b_1,\ldots,b_k)^T$ with
\[
   b_{l+1} = \psi( [ \log(F_{k,l}/ F_{k,l+1}) ] )
     = \frac12 \log (2\pi) - \log \Gamma \! \left( 1 - \frac{l}{k} \right)
       \quad \text{for} \quad l=0,\ldots,k-2
\]
using Proposition \ref{prop-fac-gamma}. The last element $b_k$ is given
by Theorem \ref{thm-asympt-fk}, Lemma \ref{lem-const-f1}, and Lemma
\ref{lem-fg-omega}:
\begin{align*}
   b_k &= \psi( [ \log(F_{1,0}(kn)) ] )
     = \frac14 \log (2\pi) + \log \MF_1 + \log \MA + \frac{5}{12} \log k \\
     &= \frac12 \log (2\pi) - \log \MA + \frac{1}{12} + \frac{5}{12} \log k.
\end{align*}
By Lemma \ref{lem-matrix-mk} we can solve the linear system directly with
\[
   x = \frac{1}{k} \, \widetilde{M_k} \cdot b .
\]
The first row yields
\[
   x_1 = \frac{1}{k} \, b_k + \frac{1}{k}
         \sum_{\nu=1}^{k-1} (k-\nu) \, b_\nu.
\]
On the other side, we have
\[
   x_1 = \psi( [ \log F_{k,0} ]) = \log \MF_k +
         \frac14 \log (2\pi) + k \log \MA.
\]
This provides
\begin{eqnarray} \label{eqn-loc-fk}
\begin{aligned}
   \log \MF_k = & -\left( k+\frac{1}{k} \right) \log \MA
     + \left( \frac{k}{4} + \frac{1}{2k} - \frac{1}{2} \right)
     \log (2\pi) \\
   & + \frac{5}{12k} \log k
     + \frac{1}{12k} - \sum_{\nu=2}^{k-1} \frac{\nu-1}{k}
     \log \Gamma \! \left( \frac{\nu}{k} \right)
\end{aligned}
\end{eqnarray}
after some rearranging of terms. By Euler's formula, see Proposition
\ref{prop-euler-gamma}, we have
\begin{equation} \label{eqn-loc-euler-prod-gamma}
   \frac{1}{k} \sum_{\nu=1}^{k-1} \log \Gamma \!
     \left( \frac{\nu}{k} \right)
     = \left( \frac12 - \frac{1}{2k} \right) \log (2\pi)
     - \frac{1}{2k} \log k.
\end{equation}
Finally, substituting \eqref{eqn-loc-euler-prod-gamma}
into \eqref{eqn-loc-fk} yields the result.
\end{proof}

\begin{remark} \label{rem-formula-fk}
Although the formula for $\MF_k$ has an elegant short form, one might also use
\eqref{eqn-loc-fk} instead, since this formula omits the value $\Gamma(1/k)$.
Thus we easily obtain the value of $\MF_2$ from \eqref{eqn-loc-fk} at once:
$\MF_2 = (2\pi)^{\frac14} \, 2^{\frac{5}{24}} \,
 e^{\frac{1}{24}} / \MA^{\frac{5}{2}}$.
\end{remark}

\begin{corl}
Asymptotically, we have
\[
   \prod_{\nu=1}^{n-1} \Gamma \! \left( \frac{\nu}{n} \right)^\nu
     \,\,\sim\,\, \frac{e^{\frac{1-\gamma}{12}}}{\MA}
     \left( \frac{(2\pi)^{\frac14}}{\MA} \right)^{n^2} \Big/ \,\,
     n^{\frac{1}{12}} \quad \text{as} \quad n \to \infty
\]
with the constants $e^{\frac{1-\gamma}{12}}/\MA = 0.8077340270...$ and
$(2\pi)^{\frac14} / \MA = 1.2345601953...$.
\end{corl}

\begin{proof}
On the one hand, we have by \eqref{eqn-log-k-fk} that
\[
   n \log \MF_n = \frac{\gamma}{12} + O ( n^{-2} ).
\]
On the other hand, Theorem \ref{thm-formula-fk} provides that
\[
   n \log \MF_n =
     - \left( n^2+1 \right) \log \MA + \frac{1}{12}
     - \frac{1}{12} \log n + \frac{n^2}{4} \log (2\pi)
     - \sum_{\nu=1}^{n-1} \nu \log \Gamma \! \left( \frac{\nu}{n} \right).
\]
Combining both formulas easily gives the result.
\end{proof}

Since we have derived exact expressions for the constants $\MF_k$, we can
improve the calculation of $\MF_\infty$. The divergent sum of $\MF_\infty$,
given in Theorem \ref{thm-asympt-fk}, is not suitable to determine a value
within a given precision, but we can use this sum in a modified way. Note
that we cannot use the limit formula
\[
   \log \MF_\infty =
     \lim_{n \to \infty} \left( \sum_{k=1}^n \log \MF_k
     - \frac{\gamma}{12} \log n \right)
\]
without a very extensive calculation, because the sequence
$\gamma_n = H_n - \log n$ converges too slowly. Moreover, the computation of
$\MF_k$ involves the computation of the values $\Gamma(\nu/k)$. This becomes
more difficult for larger $k$.

\begin{prop} \label{prop-finfty-calc}
Let $m,n$ be positive integers. Assume that $m > 2$ and the constants
$\MF_k$ are given by exact expressions for $k=1,\ldots,n$.
Define the computable values $\eta_k \in (0,1)$ implicitly by
\[
   \log \MF_k = \frac{\gamma}{12 k} +
     \sum_{j=2}^{m-1} \frac{B_{2j} \, \zeta(2j-1)}{2j (2j-1) \, k^{2j-1}}
     + \eta_k \frac{B_{2m} \, \zeta(2m-1)}{2m (2m-1) \, k^{2m-1}}.
\]
Then
\[
   \log \MF_\infty = \frac{\gamma^2}{12} +
     \sum_{j=2}^{m-1} \frac{B_{2j} \, \zeta(2j-1)^2}{2j (2j-1)}
      + \theta_{n,m} \, \frac{B_{2m} \, \zeta(2m-1)^2}{2m (2m-1)}
\]
with $\theta_{n,m} \in ( \theta_{n,m}^{\min}, \theta_{n,m}^{\max} )
\subset (0,1)$ where
\[
   \theta_{n,m}^{\min} = \zeta(2m-1)^{-1} \sum_{k=1}^{n}
     \frac{\eta_k}{k^{2m-1}},
   \qquad \theta_{n,m}^{\max} = 1 + \zeta(2m-1)^{-1}
     \sum_{k=1}^{n} \frac{\eta_k-1}{k^{2m-1}}.
\]
The error bound for the remainder of the divergent sum of
$\log \MF_\infty$ is given by
\[
   \theta_{n,m}^{\mathrm{err}}
     = \left( 1 - \zeta(2m-1)^{-1} \sum_{k=1}^{n}
     \frac{1}{k^{2m-1}} \right)
     \frac{|B_{2m}| \, \zeta(2m-1)^2}{2m (2m-1)}.
\]
\end{prop}

\begin{proof}
Let $n \geq 1 $ and $m > 2$ be fixed integers. The divergent sums for
$\log \MF_k$ and $\log \MF_\infty$ are given by Theorem \ref{thm-asympt-fk}.
Since we require exact expressions for $\MF_k$, we can compute the values
$\eta_k$ for $k=1,\ldots,n$. We define
\[
   \eta_{m,k} = \eta'_{m,k} = \eta_k \quad \text{for} \quad k=1,\ldots,n
\]
and
\[
   \eta_{m,k} = 0, \quad \eta'_{m,k} = 1 \quad \text{for} \quad k > n.
\]
We use \eqref{eqn-sum-log-fk-1} and \eqref{eqn-sum-log-fk-2} to derive
the bounds:
\[
   \theta_{n,m}^{\min} = \zeta(2m-1)^{-1} \sum_{k=1}^\infty
     \frac{\eta_{m,k}}{k^{2m-1}}
     \,\,<\,\, \theta_{n,m} \,\,<\,\,
     \zeta(2m-1)^{-1} \sum_{k=1}^\infty \frac{\eta'_{m,k}}{k^{2m-1}}
     = \theta_{n,m}^{\max}.
\]
We obtain the suggested formulas for $\theta_{n,m}^{\min}$ and
$\theta_{n,m}^{\max}$ by evaluating the sums with $\eta_{m,k} = 0$, resp.\
$\eta'_{m,k} = 1$, for $k>n$. The error bound is given by the difference of
the absolute values of the minimal and maximal remainder. Therefore
\[
   \theta_{n,m}^{\mathrm{err}}
     = (\theta_{n,m}^{\max} - \theta_{n,m}^{\min}) R
     = \left( 1 - \zeta(2m-1)^{-1} \sum_{k=1}^{n}
       \frac{1}{k^{2m-1}} \right) R
\]
with $R = |B_{2m}| \zeta(2m-1)^2/ 2m (2m-1)$.
\end{proof}

\begin{result}
Exact expressions for $\MF_k$:
\[
\begin{aligned}
   \MF_1 &= (2\pi)^{\frac14} \, e^{\frac{1}{12}} / \MA^2, &
   \MF_2 &= (2\pi)^{\frac14} \, 2^{\frac{5}{24}} \,
     e^{\frac{1}{24}} / \MA^{\frac{5}{2}}, \\
   \MF_3 &= (2\pi)^{\frac{5}{12}} \, 3^{\frac{5}{36}} \,
     e^{\frac{1}{36}} / \MA^{\frac{10}{3}} \,
     \Gamma\! \left( \textstyle\frac{2}{3} \right)^\frac13, &
   \MF_4 &= (2\pi)^{\frac{1}{2}} \, 2^{\frac{1}{3}} \,
     e^{\frac{1}{48}} / \MA^{\frac{17}{4}} \,
     \Gamma\! \left( \textstyle\frac{3}{4} \right)^\frac12.
\end{aligned}
\]

We have computed the constants $\MF_k$ by their exact expression.
Moreover, we have determined the index $m$ of the smallest remainder
of their asymptotic divergent series and the resulting error bound
given by Theorem \ref{thm-asympt-fk}.

\begin{center}
\begin{tabular}{|c|l|r|l|} \hline
  Constant & \multicolumn{1}{c|}{Value} & $m$ &
    \multicolumn{1}{c|}{Error bound} \\\hline \hline
  $\MF_1$ & 1.04633506677050318098... &  4 & $6.000\cdot10^{-4}$ \\
  $\MF_2$ & 1.02393741163711840157... &  7 & $7.826\cdot10^{-7}$ \\
  $\MF_3$ & 1.01604053706462099128... & 10 & $1.198\cdot10^{-9}$ \\
  $\MF_4$ & 1.01204589802394464624... & 13 & $1.948\cdot10^{-12}$\\
  $\MF_5$ & 1.00963997283647705086... & 16 & $3.272\cdot10^{-15}$\\
  $\MF_6$ & 1.00803362724207326544... & 20 & $5.552\cdot10^{-18}$\\\hline
\end{tabular}
\end{center}
\smallskip

The weak interval of $\MF_\infty$ is given by Theorem \ref{thm-asympt-fk}. The
second value is derived by Proposition \ref{prop-finfty-calc} with parameters
$m=17$ and $n=7$. Thus, exact expressions of $\MF_1,\ldots,\MF_7$ are needed to
compute $\MF_\infty$ within the given precision.
\smallskip

\begin{center}
\begin{tabular}{|c|l|r|l|} \hline
  Constant & \multicolumn{1}{c|}{Value / Interval} & $m$ &
    \multicolumn{1}{c|}{Error bound} \\\hline \hline
  $\MF_\infty$ & \multicolumn{1}{c|}{(1.02428, 1.02491)} & 4 &
    $6.050\cdot10^{-4}$ \\
  $\MF_\infty$ & 1.02460688265559721480... & 17 &
    $6.321\cdot10^{-22}$ \\\hline
\end{tabular}
\end{center}
\end{result}

\section{Products of Bernoulli numbers}

Using results of the previous sections, we are now able to consider several
products of Bernoulli numbers and to derive their asymptotic expansions and
constants.

\begin{theorem} \label{thm-asympt-bn}
Asymptotically, we have
\begin{eqnarray*}
   \prod_{\nu=1}^n |B_{2\nu}| &\sim& \MB_1
     \left( \frac{n}{\pi e^{3/2}} \right)^{n(n+1)} \,
     (16 \pi n)^{\frac{n}{2}} \, n^{\frac{11}{24}}
     \quad \text{as} \quad n \to \infty, \\
   \prod_{\nu=1}^n \frac{|B_{2\nu}|}{2\nu} &\sim& \MB_2
     \left( \frac{n}{\pi e^{3/2}} \right)^{n^2}
     \left( \frac{4n}{\pi e} \right)^{\frac{n}{2}} \Big/ \,
     n^{\frac{1}{24}} \quad \text{as} \quad n \to \infty
\end{eqnarray*}
with the constants
\[
   \begin{array}{c@{\,\,=\,\,}l@{\,\,=\,\,}l}
     \MB_1 & \MC_2 \MF_2 \MA^2 (2\pi)^{\frac14}
       & \MC_2 \, (2\pi)^{\frac12} \, 2^{\frac{5}{24}} \,
         e^{\frac{1}{24}} / \MA^{\frac{1}{2}}, \\
     \MB_2 & \MC_2 \MF_2 \MA^2 / (2\pi)^{\frac14}
       & \MC_2 \, 2^{\frac{5}{24}} \,
         e^{\frac{1}{24}} / \MA^{\frac{1}{2}}.
   \end{array}
\]
\end{theorem}

\begin{proof}
By Euler's formula \eqref{eqn-euler-zeta} for $\zeta(2\nu)$ and Lemma
\ref{lem-const-c1} we obtain
\[
   \prod_{\nu=1}^n |B_{2\nu}| \,\,\sim\,\,
     \MC_2 \, \prod_{\nu=1}^n \frac{2 \cdot (2\nu)!}{(2\pi)^{2\nu}}
     \,\,\sim\,\, \MC_2 \, 2^n (2\pi)^{-n(n+1)} \prod_{\nu=1}^n (2\nu)!
     \quad \text{as} \quad n \to \infty.
\]
Theorem \ref{thm-asympt-fk} states for $k=2$ that
\[
   \prod_{\nu=1}^n (2 \nu)!  \,\,\sim\,\,
     \MF_2 \, \MA^2 \, (2\pi)^{\frac14}
     \left( \frac{2n}{e^{3/2}} \right)^{n(n+1)}
     \left( 4\pi n \right)^{\frac{n}{2}}
     n^{\frac{11}{24}} \quad \text{as} \quad n \to \infty.
\]
The expression for $\MF_2$ is given in Remark \ref{rem-formula-fk}.
Combining both asymptotic formulas above gives the first suggested formula.
It remains to evaluate the following product:
\[
   \prod_{\nu=1}^n (2\nu) = 2^n \, n! \,\,\sim\,\, (2\pi)^{\frac12}
     \left( \frac{2n}{e} \right)^n n^{\frac12}
     \quad \text{as} \quad n \to \infty.
\]
After some rearranging of terms we then obtain the second suggested formula.
\end{proof}

\begin{remark} \label{rem-milnor}
Milnor and Husemoller \cite[pp.~49--50]{Milnor73} give the following
asymptotic formula without proof:
\begin{equation} \label{eqn-mil-1}
   \prod_{\nu=1}^n |B_{2\nu}| \,\,\sim\,\, \MB' \,
     n! \, 2^{n+1} \, F(2n+1) \quad \text{as} \quad n \to \infty
\end{equation}
where
\begin{equation} \label{eqn-mil-2}
   F(n) = \left( \frac{n}{2\pi e^{3/2}} \right)^{\frac{n^2}{4}}
     \left( \frac{8\pi e}{n} \right)^{\frac{n}{4}} \Big/ \, n^{\frac{1}{24}}
\end{equation}
and $\MB' \approx 0.705$ is a certain constant. This constant is related to
the constant $\MB_2$.
\end{remark}

\begin{prop}
The constant $\MB'$ is given by
\[
   \MB' = 2^{\frac{1}{24}} \, 2^{-\frac{3}{2}} \, \MB_2
     = \MC_2 \, e^{\frac{1}{24}} / 2^{\frac54} \MA^{\frac12}
     = 0.7048648734... \,.
\]
\end{prop}

\begin{proof}
By Theorem \ref{thm-asympt-bn} we have
\begin{equation} \label{eqn-asym-gn}
   \prod_{\nu=1}^n \frac{|B_{2\nu}|}{2\nu}
     \,\,\sim\,\, \MB_2 \, G(n) \quad \text{as} \quad n \to \infty
\end{equation}
with
\[
   G(n) = \left( \frac{n}{\pi e^{3/2}} \right)^{n^2}
     \left( \frac{4n}{\pi e} \right)^{\frac{n}{2}}
     \Big/ \, n^{\frac{1}{24}}.
\]
We observe that \eqref{eqn-mil-1} and \eqref{eqn-asym-gn} are equivalent
so that
\[
   2 \, \MB' F(2n+1) \,\,\sim\,\, \MB_2 \, G(n)
   \quad \text{as} \quad n \to \infty.
\]
We rewrite \eqref{eqn-mil-2} in the suitable form
\[
   F(2n+1) = \left( \frac{n+\frac12}{\pi e^{3/2}} \right)^{n^2+n+\frac14}
     \left( \frac{4\pi e}{n+\frac12} \right)^{\frac{n}{2}+\frac14} \Big/ \,
     2^{\frac{1}{24}} \, \left( n+\frac12 \right)^{\frac{1}{24}}.
\]
Hence, we easily deduce that
\[
   G(n)/F(2n+1) = \left( 1 + \frac{1}{2n}
     \right)^{-n^2-\frac{n}{2}+\frac{1}{24}} \,
     e^{\frac{n}{2}} \, 2^{\frac{1}{24}} \,
     \left( \frac{e^{1/2}}{4} \right)^{\frac14}.
\]
It is well known that
\[
   \lim_{n \to \infty} \left( 1 + \frac{x}{n} \right)^n = e^x
     \qquad \text{and} \qquad
     \lim_{n \to \infty} e^{-xn} \left( 1 + \frac{x}{n} \right)^{n^2}
     = e^{-\frac{x^2}{2}}.
\]
Evaluating the asymptotic terms, we get
\[
   2 \, \MB' / \MB_2 \,\,\sim\,\, G(n) / F(2n+1) \,\,\sim\,\,
     e^{\frac18} \, e^{-\frac14} \, 2^{\frac{1}{24}} \,
     e^{\frac18} \, 2^{-\frac12} \quad \text{as} \quad n \to \infty,
\]
which finally yields $\MB' = 2^{\frac{1}{24}} \, 2^{-\frac{3}{2}} \, \MB_2$.
\end{proof}

\begin{theorem}
The Minkowski--Siegel mass formula asymptotically states for positive
integers $n$ with $4 \pdiv n$ that
\[
   M(2n) = \frac{|B_n|}{2n} \,\, \prod_{\nu=1}^{n-1} \frac{|B_{2\nu}|}{4\nu}
     \,\,\sim\,\, \MB_3
     \left( \frac{n}{\pi e^{3/2}} \right)^{n^2} \Big/ \,
       \left( \frac{4n}{\pi e} \right)^{\frac{n}{2}} \, n^{\frac{1}{24}}
       \quad \text{as} \quad n \to \infty
\]
with $\MB_3 = \sqrt{2} \, \MB_2$.
\end{theorem}

\begin{proof}
Let $n$ always be even. By Proposition \ref{prop-gamma-stirling} and
\eqref{eqn-euler-zeta} we have
\[
   2^{-n} \left| \frac{B_n/n}{B_{2n}/2n} \right|
     = 2 \frac{\zeta(n)}{\zeta(2n)}
     \pi^n \frac{n!}{(2n)!} \,\,\sim\,\,
     \sqrt{2} \left( \frac{4n}{\pi e} \right)^{-n}
     \quad \text{as} \quad n \to \infty,
\]
since $\zeta(n)/\zeta(2n) \sim 1$ and
\[
   \log \left( \frac{n!}{(2n)!} \right)
     \,\,\sim\,\, n - n \log n - \left( 2n+\frac12 \right) \log 2
     \quad \text{as} \quad n \to \infty .
\]
We finally use Theorem \ref{thm-asympt-bn} and \eqref{eqn-asym-gn} to obtain
\[
   M(2n) = 2^{-n} \left| \frac{B_n/n}{B_{2n}/2n} \right|
     \prod_{\nu=1}^{n} \frac{|B_{2\nu}|}{2\nu} \,\,\sim\,\,
     \sqrt{2} \, \MB_2 \left( \frac{4n}{\pi e} \right)^{-n} G(n)
     \quad \text{as} \quad n \to \infty,
\]
which gives the result.
\end{proof}

\begin{result}
The constants $\MB'$, $\MB_\nu$ $(\nu = 1,2,3)$ mainly depend on the constant
$\MC_2$ and the Glaisher--Kinkelin constant $\MA$.
\smallskip

\begin{center}
\begin{tabular}{|c|c|c|} \hline
  Constant & Expression & Value \\\hline \hline
  $\MA$    & & 1.28242712910062263687... \\
  $\MC_2$  & & 1.82101745149929239040... \\
  $\MB_1$  & $\MC_2 (2\pi)^{\frac12} \, 2^{\frac{5}{24}} \,
     e^{\frac{1}{24}} / \MA^{\frac{1}{2}}$
     & 4.85509664652226751252... \\
  $\MB_2$  & $\MC_2 2^{\frac{5}{24}} \,
     e^{\frac{1}{24}} / \MA^{\frac{1}{2}}$
     & 1.93690332773294192068... \\
  $\MB_3$  & $\MC_2 2^{\frac{17}{24}} \,
     e^{\frac{1}{24}} / \MA^{\frac{1}{2}}$
     & 2.73919495508550621998... \\
  $\MB'$   & $\MC_2 \, e^{\frac{1}{24}} / 2^{\frac54} \MA^{\frac12}$
     & 0.70486487346802031057... \\ \hline
\end{tabular}
\end{center}
\end{result}

\section{Generalizations}

In this section we derive a generalization of Theorem \ref{thm-asympt-fk}.
The results show the structure of the constants $\MF_k$ and the generalized
constants $\MF_{r,k}$, which we shall define later, in a wider context. For
simplification we introduce the following definitions which arise from the
Euler-Maclaurin summation formula.

The sum of consecutive integer powers is given by the well known formula
\[
   \sum_{\nu=0}^{n-1} \nu^r = \frac{B_{r+1}(n) - B_{r+1}}{r+1}
     = \sum_{j=0}^r \binom{r}{j} B_{r-j} \frac{n^{j+1}}{j+1},
       \quad r \geq 0,
\]
where $B_{m}(x)$ is the $m$th Bernoulli polynomial. Now, the Bernoulli number
$B_1=-\frac12$ is responsible for omitting the last power $n^r$ in the
summation above. Because we further need the summation up to $n^r$, we change
the sign of $B_1$ in the sum as follows:
\[
   S_r(n) = \sum_{\nu=1}^{n} \nu^r
     = \sum_{j=0}^r \binom{r}{j} (-1)^{r-j} B_{r-j} \frac{n^{j+1}}{j+1},
       \quad r \geq 0 .
\]
This modification also coincides with
\[
   \zeta(-n) = (-1)^{n+1} \frac{B_{n+1}}{n+1}
\]
for nonnegative integers $n$. We define the extended sum
\[
   S_r(n;f(\diamond)) = \sum_{j=0}^r \binom{r}{j}
     (-1)^{r-j} B_{r-j} \frac{n^{j+1}f(j+1)}{j+1},
       \quad r \geq 0,
\]
where the symbol $\diamond$ is replaced by the index $j+1$ in the sum.
Note that $S_r$ is linear in the second parameter, i.e.,
\[
   S_r(n;\alpha+\beta f(\diamond) )
     = \alpha S_r(n) + \beta S_r(n; f(\diamond) ).
\]
Finally we define
\[
   D_k(x) = \sum_{j \geq 1} \divsum \BN_{2j,k} \, x^{-(2j-1)}
   \quad \text{where} \quad
   \BN_{m,k} = \frac{B_m}{m(m-1)k^{m-1}} .
\]

\begin{theorem} \label{thm-gen-glaisher}
Let $r$ be a nonnegative integer. Then
\[
   \prod_{\nu=1}^{n} \nu^{\, \nu^r} \,\,\sim\,\, \MA_r \, Q_r(n)
   \quad \text{as} \quad n \to \infty,
\]
where $\MA_r$ is the generalized Glaisher--Kinkelin constant defined by
\[
   \log \MA_r = -\zeta(-r) \, H_r - \zeta'(-r) .
\]
Moreover, $\log Q_r \in \Omega_{r+1}$ with
\[
   \log Q_r(n) = ( S_r(n) - \zeta(-r) ) \log n
     + S_r(n;H_r-H_{\diamond}) .
\]
\end{theorem}

\begin{proof}
This formula and the constants easily follow from a more general formula for
real $r > -1$ given in \cite[9.28, p.~595]{Graham94} and after some rearranging
of terms.
\end{proof}

\begin{remark}
The case $r=0$ reduces to Stirling's approximation of $n!$ with
$\MA_0 = \sqrt{2\pi}$. The case $r=1$ gives the usual Glaisher--Kinkelin
constant $\MA_1 = \MA$. The expression $S_r(n;H_r-H_{\diamond})$ does not
depend on the definition of $B_1$, since the term with $B_1$ is cancelled
in the sum. Graham, Knuth, and Patashnik \cite[9.28, p.~595]{Graham94} notice
that the constant $- \zeta'(-r)$ has been determined in a book of de Bruijn
\cite[\S 3.7]{deBruijn70} in 1970. The theorem above has a long history.
In 1894 Alexeiewsky \cite{Alexeiewsky94} gave the identity
\[
   \prod_{\nu=1}^{n} \nu^{\, \nu^r}
     = \exp\left( \zeta'(-r,n+1) - \zeta'(-r) \right)
\]
where $\zeta'(s,a)$ is the partial derivative of the Hurwitz zeta function with
respect to the first variable. Between 1903 and 1913, Ramanujan recorded in his
notebooks \cite[Entry~27, pp.~273--276]{Berndt85} (the first part was published
and edited by Berndt \cite{Berndt85} in 1985) an asymptotic expansion for real
$r > -1$ and an analytic expression for the constant $C_r = - \zeta'(-r)$.
However, Ramanujan only derived closed expressions for $C_0$ and
$C_{2r}$ ($r \geq 1$) in terms of $\zeta(2r+1)$; see \eqref{eqn-zeta-deriv-even}
below. In 1933 Bendersky \cite{Bendersky33} showed that there exist certain
constants $\MA_r$. Since 1980, several others have investigated the asymptotic
formula, including MacLeod \cite{MacLeod82}, Choudhury \cite{Choudhury95}, and
Adamchik \cite{Adamchik98, Adamchik01}.
\end{remark}

\begin{theorem} \label{thm-gen-asympt-fk}
Let $k, r$ be integers with $k \geq 1$ and $r \geq 0$. Then
\[
   \prod_{\nu=1}^n (k \nu)!^{\,\nu^r}  \,\,\sim\,\,
     \MF_{r,k} \, \MA_r^{\frac12} \, \MA_{r+1}^k \,
     P_{r,k}(n) \, Q_r(n)^{\frac12} \, Q_{r+1}(n)^k
     \quad \text{as} \quad n \to \infty.
\]
The constants $\MF_{r,k}$ and functions $P_{r,k}$ satisfy that $\lim\limits_{k
\to \infty} \MF_{r,k} = 1$ and $\log P_{r,k} \in \Omega_{r+2}$ where
\begin{align*}
   \log P_{r,k}(n) = & \, \frac12 S_r(n) \log(2\pi k)
     + k\, S_{r+1}(n) \log (k/e) \\
   & + \BN_{r+2,k} \log n + \sum_{j=1}^{\lfloor \frac{r+1}{2} \rfloor}
     \BN_{2j,k} \, S_{r+1-2j}(n).
\end{align*}
The constants $\MA_r$ and functions $Q_r$ are defined as in
Theorem \ref{thm-gen-glaisher}.
\end{theorem}

The determination of exact expressions for the constants $\MF_{r,k}$
seems to be a very complicated and extensive task in the case $r > 0$.
The next theorem gives a partial result for $k=1$ and $r \geq 0$.

\begin{theorem} \label{thm-gen-fr1}
Let $r$ be a nonnegative integer. Then
\[
   \log \MF_{r,1} = \frac12 \log \MA_r - \log \MA_{r+1}
     + S_r(1;\BN_{1+\diamond,1}-\log \MA_\diamond).
\]
Case $r=0$:
\[
   \log \MF_{r,1} = \frac{1}{12} + \frac12 \log \MA_0 - 2 \log \MA_1.
\]
Case $r>0$:
\[
   \log \MF_{r,1} = \alpha_{r,0} +
     \sum_{j=1}^{r+1} \alpha_{r,j} \, \log \MA_j
\]
where
\[
   \alpha_{r,j} = \left\{
     \begin{array}{rll}
       \frac{B_{r+1}}{2r(r+1)},
         & r \not\equiv j \ \mods{2}, & j = 0 ; \\
       \sum\limits_{j=0}^{r} \! \binom{r}{j} \frac{B_{r-j} \,
         B_{j+2}}{(j+1)^2 (j+2)},
         & r \equiv j \ \mods{2}, & j = 0 ; \\
       -\delta_{r+1,j}- \binom{r+1}{j} \frac{B_{r+1-j}}{r+1},
         & r \not\equiv j \ \mods{2}, & j > 0 ; \vspace*{1ex} \\
       0, & r \equiv j \ \mods{2}, & j > 0
     \end{array} \right.
\]
and $\delta_{i,j}$ is Kronecker's delta.
\end{theorem}

\begin{proof}[Proof of Theorem~\ref{thm-gen-asympt-fk}]
Let $k$ and $r$ be fixed.
We extend the proof of Theorem \ref{thm-asympt-fk}.
From \eqref{eqn-log-kn-fac} we have
\begin{equation} \label{eqn-log-kn-fac-dk}
   \log (k\nu)! = \frac12 \log (2\pi k)
     + k \nu \log \left( \frac{k}{e} \right)
     + \left( k \nu + \frac12 \right) \log \nu + D_k(\nu).
\end{equation}
The summation yields
\[
   \sum_{\nu=1}^n \nu^r \log (k\nu)! = F_1(n) + F_2(n) + F_3(n)
\]
where
\begin{align*}
   F_1(n) &= \frac{1}{2} S_r(n) \log (2\pi k) + k S_{r+1}(n) \log (k/e), \\
   F_2(n) &= k \sum_{\nu=1}^n \nu^{r+1} \log \nu + \frac12
             \sum_{\nu=1}^n \nu^r \log \nu, \\
   F_3(n) &= \sum_{\nu=1}^n \nu^r D_k(\nu).
\end{align*}
Theorem \ref{thm-gen-glaisher} provides
\[
   F_2(n) = k \left( \log \MA_{r+1} + \log Q_{r+1}(n) \right) +
     \frac12 \left( \log \MA_{r} + \log Q_{r}(n) \right) + O(n^{-\delta})
\]
with some $\delta > 0$. Let $R=\lfloor \frac{r+1}{2} \rfloor$.
By definition we have
\[
   x^r D_k(x) = \sum_{j=1}^{R} \BN_{2j,k} \, x^{r+1-2j}
     + \sum_{j > R} \divsum \BN_{2j,k} \, x^{r+1-2j}
   =: E_1(x) + E_2(x).
\]
Therewith we obtain that
\[
   F_3(n) = \sum_{j=1}^{R} \BN_{2j,k} \, S_{r+1-2j}(n)
     + \sum_{\nu=1}^n E_2(\nu).
\]
For the second sum above we consider two cases. We use similar arguments
which we have applied to \eqref{eqn-log-fk-1} and \eqref{eqn-log-fk-2}.
If $r$ is odd, then
\begin{equation} \label{eqn-e2-r-odd}
   \lim_{n \to \infty} \sum_{\nu=1}^n E_2(\nu)
     = \sum_{j > R} \divsum \BN_{2j,k} \, \zeta(2j-(r+1)).
\end{equation}
Note that $\BN_{r+2,k}=0$ in that case.
If $r$ is even, then we have to take care of the term $\nu^{-1}$.
This gives
\begin{equation} \label{eqn-e2-r-even}
   \lim_{n \to \infty} \left( \sum_{\nu=1}^n E_2(\nu)
     - \BN_{r+2,k} \log n \right)
   = \gamma \, \BN_{r+2,k} + \sum_{j > R+1} \!\!\!\divsum
     \BN_{2j,k} \, \zeta(2j-(r+1)).
\end{equation}
The right hand side of \eqref{eqn-e2-r-odd}, resp.\ \eqref{eqn-e2-r-even},
defines the constant $\log \MF_{r,k}$. Finally we have to collect all results
for $F_1$, $F_2$, and $F_3$. This gives the constants and the function
$P_{r,k}$. It remains to show that $\lim_{k \to \infty} \log \MF_{r,k} = 0$.
This follows by $\BN_{2j,k} \to 0$ as $k \to \infty$.
\end{proof}

The following lemma gives a generalization of
Equation \eqref{eqn-prod-fac-f1} in Lemma \ref{lem-const-f1}.
After that we can give a proof of Theorem~\ref{thm-gen-fr1}.

\begin{lemma}
Let $n, r$ be integers with $n \geq 1$ and $r \geq 0$. Then
\begin{equation} \label{eqn-prod-fac-gen}
   n!^{S_r(n)} \prod_{\nu=1}^n \nu^{\,\nu^r}
     = \prod_{\nu=1}^n \nu!^{\,\nu^r} \, \prod_{\nu=1}^n \nu^{S_r(\nu)}.
\end{equation}
\end{lemma}

\begin{proof}
We regard the following enumeration scheme which can be easily extended
to $n$ rows and $n$ columns:
\[
   \begin{array}{ccc}
     1^{1^r} & \framebox{$2^{1^r}\!$} & \framebox{$3^{1^r}\!$} \\
     1^{2^r} & 2^{2^r} & \framebox{$3^{2^r}\!$} \\
     1^{3^r} & 2^{3^r} & 3^{3^r} \\
   \end{array}
\]
The product of all elements, resp.\ non-framed elements, in the $\nu$th row
equals $n!^{\,\nu^r}$, resp.\ $\nu!^{\,\nu^r}$. The product of the framed
elements in the $\nu$th column equals $\nu^{S_r(\nu-1)}$. Thus
\[
   n!^{S_r(n)} = \prod_{\nu=1}^n \nu!^{\,\nu^r} \,
     \prod_{\nu=1}^n \nu^{S_r(\nu)-\nu^r}. \qedhere
\]
\end{proof}

\begin{proof}[Proof of Theorem~\ref{thm-gen-fr1}]
Let $r \geq 0$.
We take the logarithm of \eqref{eqn-prod-fac-gen} to obtain
\begin{equation} \label{eqn-gen-fr1-0}
   F_1(n) + F_2(n) = F_3(n) + F_4(n)
\end{equation}
where
\[
\begin{aligned}
   F_1(n) &= S_r(n) \log n!, &
   F_2(n) &= \sum_{\nu=1}^n \nu^r \log \nu, \\
   F_3(n) &= \sum_{\nu=1}^n \nu^r \log \nu!, &
   F_4(n) &= \sum_{\nu=1}^n S_r(\nu) \log \nu.
\end{aligned}
\]
Next we consider the asymptotic expansions $\tilde{F}_j$ of
the functions $F_j$ ($j=1,\ldots,4$) when $n \to \infty$.
We further reduce the functions $\tilde{F}_j$ via the maps
\[
   C(\RR^+;\RR) \, \stackrel{[\ ]}{\longrightarrow} \,
     \Omega_\infty \, \stackrel{\psi}\longrightarrow \, \RR
\]
to the constant terms which are the asymptotic constants of
$[\tilde{F}_j]$ in $\Omega_\infty$.
Consequently \eqref{eqn-gen-fr1-0} turns into
\begin{equation} \label{eqn-gen-fr1-1}
   \psi([\tilde{F}_1]) + \psi([\tilde{F}_2])
     = \psi([\tilde{F}_3]) + \psi([\tilde{F}_4]).
\end{equation}
We know from Theorem \ref{thm-gen-glaisher} and
Theorem \ref{thm-gen-asympt-fk} that
\[
   \psi([\tilde{F}_2]) = \log \MA_r \quad \text{and} \quad
     \psi([\tilde{F}_3]) = \log \MF_{r,1} + \frac12 \log \MA_r
     + \log \MA_{r+1}.
\]
For $\tilde{F}_4$ we derive the expression
\begin{equation} \label{eqn-gen-fr1-2}
   \psi([\tilde{F}_4]) = S_r(1;\log \MA_\diamond ),
\end{equation}
since each term $s_j \nu^j$ in $S_r(\nu)$ produces the term $s_j \log \MA_j$.
It remains to evaluate $\tilde{F}_1$.
According to \eqref{eqn-log-kn-fac-dk} we have
\[
   \log n! = \frac12 \log (2\pi) - n + \left( n+\frac12 \right)
     \log n + D_1(n) =: E(n) + D_1(n).
\]
Thus
\[
   \tilde{F}_1(x) = S_r(x) E(x) + S_r(x) D_1(x).
\]
Since $S_r E \in \Omega_\infty$ has no constant term, we deduce that
\[
   \psi([\tilde{F}_1]) = \psi([S_r D_1]) = S_r(1;\BN_{1+\diamond,1}).
\]
The latter equation is similarly derived as \eqref{eqn-gen-fr1-2}, whereas
we regard the constant terms of the product of the polynomial $S_r$ and the
Laurent series $D_1$. From \eqref{eqn-gen-fr1-1} we finally obtain
\[
   \log \MF_{r,1} = \frac12 \log \MA_r - \log \MA_{r+1}
     + S_r(1;\BN_{1+\diamond,1}-\log \MA_\diamond).
\]
Now, we shall evaluate the expression above. For $r=0$ we get
\[
    \log \MF_{0,1} = \frac{1}{12} + \frac12 \log \MA_0 - 2 \log \MA_1,
\]
since
\[
   S_0(1;\BN_{1+\diamond,1}-\log \MA_\diamond) = \BN_{2,1}-\log \MA_1
     = \frac{1}{12}-\log \MA_1.
\]
For now, let $r > 0$. We may represent $\log \MF_{r,1}$ in terms of
$\log \MA_j$ as follows:
\[
   \log \MF_{r,1} = \alpha_{r,0} +
     \sum_{j=1}^{r+1} \alpha_{r,j} \, \log \MA_j.
\]
The term $\alpha_{r,0}$ is given by
\[
   \alpha_{r,0} = S_r(1;\BN_{1+\diamond,1})
     = \sum_{j=0}^r \binom{r}{j} (-1)^{r-j} B_{r-j} \frac{\BN_{j+2,1}}{j+1}
\]
where the sum runs over even $j$, since $\BN_{j+2,1} = 0$ for odd $j$.
If $r$ is odd, then the sum simplifies to the term $B_{r+1}/2r(r+1)$.
Otherwise we derive for even $r$ that
\[
   \alpha_{r,0} = \sum_{j=0}^{r} \! \binom{r}{j} \frac{B_{r-j} \,
     B_{j+2}}{(j+1)^2 (j+2)}.
\]
It remains to determine the coefficients $\alpha_{r,j}$ for $r+1 \geq j \geq 1$.
Since $\frac12 x^r - x^{r+1} - S_r(x)$ is an odd, resp.\ even, polynomial
for even, resp.\ odd, $r > 0$, this property transfers in a similar way to
$\frac12 \log \MA_r - \log \MA_{r+1} - S_r(1;\log \MA_\diamond)$, such that
$\alpha_{r,j} = 0$ when $2 \pdiv r-j$. Otherwise we get
\begin{equation} \label{eqn-gen-fr1-3}
   \alpha_{r,j} = - \binom{r}{j-1} \frac{B_{r-(j-1)}}{j} - \delta_{r+1,j}
     = - \binom{r+1}{j} \frac{B_{r+1-j}}{r+1} - \delta_{r+1,j}
\end{equation}
for $2 \notdiv r-j$,
where the term $- \log \MA_{r+1}$ is represented by $- \delta_{r+1,j}$.
\end{proof}

\begin{corl} \label{corl-fr1-odd}
Let $r$ be an odd positive integer. Then
\begin{align*}
   \log \MF_{r,1} &= - \frac{r!}{(2\pi i)^{r+1}}
     \left( \frac{\zeta(r+1)}{r}
     + \sum_{j=1}^{\frac{r-1}{2}} \zeta(r+1-2j) \zeta(2j+1)
     - \frac{(r+2)\zeta(r+2)}{2} \right) \\
   &= (-1)^{\frac{r-1}{2}} \frac{r!}{2} \left(
     \frac{|B_{r+1}|}{r(r+1)!} + \sum_{j=1}^{\frac{r-1}{2}}
     \frac{|B_{r+1-2j}|\,\zeta(2j+1)}{(r+1-2j)!\,(2\pi)^{2j}}
     - \frac{(r+2)\zeta(r+2)}{(2\pi)^{r+1}} \right).
\end{align*}
\end{corl}

\begin{proof}
As a consequence of the functional equation of $\zeta(s)$ and its derivative,
we have for even positive integers $n$, cf.\ \cite[p.~276]{Berndt85}, that
\begin{equation} \label{eqn-zeta-deriv-even}
   \log \MA_n = - \zeta'(-n) = - \frac12 \frac{n!}{(2\pi i)^n} \zeta(n+1)
\end{equation}
where the left hand side of \eqref{eqn-zeta-deriv-even} follows by definition.
Theorem \ref{thm-gen-fr1} provides
\[
   \log \MF_{r,1} = \frac{B_{r+1}}{2r(r+1)} +
     \sum_{j=1}^{\frac{r+1}{2}} \alpha_{r,2j} \, \log \MA_{2j} .
\]
Combining \eqref{eqn-gen-fr1-3} and \eqref{eqn-zeta-deriv-even} gives
the second equation above.
By Euler's formula \eqref{eqn-euler-zeta} we finally derive the first equation.
\end{proof}

\begin{remark}
For the sake of completeness, we give an analogue of
\eqref{eqn-zeta-deriv-even} for odd integers.
From the logarithmic derivatives of $\Gamma(s)$ and
the functional equation of $\zeta(s)$,
see \cite[pp.~183, 276]{Berndt85},
it follows for even positive integers $n$, that
\[
   \log \MA_{n-1} = \frac{B_n}{n} H_{n-1} - \zeta'(1-n)
     = \frac{B_n}{n}(\gamma + \log(2\pi))
       + 2 \frac{(n-1)!}{(2\pi i)^n} \zeta'(n)
\]
where
\[
   \zeta'(n) = - \sum_{\nu=2}^\infty \log(\nu) \, \nu^{-n}.
\]
However, \textsc{Mathematica} is able to compute values of $\zeta'$
for positive and negative argument values to any given precision.
\end{remark}

\begin{result}
Exact expressions for $\MF_{r,1}$ in terms of $\MA_j$:

\begin{center}
\begin{tabular}{|c|c|c|} \hline
  Constant & Expression & Value \\\hline \hline
  $\MF_{0,1}$ & $e^\frac{1}{12} \MA_0^\frac12 \MA_1^{-2}$
    & 1.04633506677050318098... \\
  $\MF_{1,1}$ & $e^\frac{1}{24} \MA_2^{-\frac32}$
    & 0.99600199446870605433... \\
  $\MF_{2,1}$ & $e^\frac{7}{540} \MA_1^{-\frac16} \MA_3^{-\frac43}$
    & 0.99904614418135586848... \\
  $\MF_{3,1}$ & $e^{-\frac{1}{720}} \MA_2^{-\frac14} \MA_4^{-\frac54}$
    & 1.00097924030236153773... \\
  $\MF_{4,1}$ & $e^{-\frac{67}{18900}} \MA_1^{\frac{1}{30}}
    \MA_3^{-\frac13} \MA_5^{-\frac65}$
    & 1.00007169725554110099... \\
  $\MF_{5,1}$ & $e^{\frac{1}{2520}} \MA_2^{\frac{1}{12}}
    \MA_4^{-\frac{5}{12}} \MA_6^{-\frac76}$
    & 0.99937792615674804266... \\\hline
\end{tabular} \smallskip
\end{center}
Exact expressions for $\MF_{r,1}$ in terms of $\zeta(2j+1)$:
\begin{align*}
   \MF_{1,1} &= \exp \left( \frac{1}{24} - \frac{3\zeta(3)}{8\pi^2}
                \right) \!, \\
   \MF_{3,1} &= \exp \left( -\frac{1}{720} - \frac{\zeta(3)}{16\pi^2}
                + \frac{15\zeta(5)}{16\pi^4} \right) \!, \\
   \MF_{5,1} &= \exp \left( \frac{1}{2520} + \frac{\zeta(3)}{48\pi^2}
                + \frac{5\zeta(5)}{16\pi^4} - \frac{105\zeta(7)}{16\pi^6}
                \right) \!.
\end{align*}
For the first 15 constants $\MF_{r,1}$ ($r=0,\ldots,14$) we find that
\[
   \max_{0 \leq r \leq 14} | \MF_{r,1} - 1 | < 0.05 ,
\]
but, e.g., $\MF_{19,1} \approx 371.61$ and
$\MF_{20,1} \approx 1.16 \cdot 10^{-7}$.
\end{result}

\section*{Acknowledgements}

The author would like to thank Steven Finch for informing about the formula of
Milnor--Husemoller and the problem of finding suitable constants; also for
giving several references. The author is also grateful to the referee for
several remarks and suggestions.

\bibliographystyle{amsplain}

\end{document}